\documentclass{article}
\usepackage{amsmath,amssymb,indentfirst}
\usepackage[latin1]{inputenc}
\usepackage{url}
\newtheorem{theorem}{Theorem}
\newtheorem{lemma}{Lemma}

\newcommand{\qed}{\hfill\mbox{\raggedright $\Box$}\medskip}

\begin{document}

\title{THE IMPLICIT FUNCTION THEOREM WHEN THE 
MATRIX  $\frac{\partial F}{\partial y}(x,y)$ IS ONLY CONTINUOUS AT THE BASE POINT}

\date{}
\maketitle

\begin{abstract} This article presents an elementary proof of the Implicit Function Theorem for differentiable maps $F(x,y)$, defined on a finite-dimensional Euclidean space, with $\frac{\partial F}{\partial y}(x,y)$ only continuous at the base point. In the case of a single scalar equation this continuity hypothesis is not required. The Inverse Function Theorem is also shown. The proofs rely on the mean-value and the intermediate-value theorems and Darboux's property (the intermediate-value property for derivatives). These proofs  avoid compactness arguments, fixed-point theorems, and integration theory. 
\end{abstract}

\vspace{0,2 cm}

\hspace{- 0,6 cm}{\sl Mathematics Subject Classification: 26B10, 
26B12}

\hspace{- 0,6 cm}{\sl Key words and phrases:} Implicit Function Theorems, Jacobians, Transformations with Several Variables, Calculus of Vector Functions.

\vspace{0,3 cm}

\section{Introduction.} The aim of this article is to present a very elementary proof of a quite well-known and generally easy to apply Implicit Function Theorem. We prove this theorem for differentiable maps $F(x,y)$ defined on a finite-dimensional Euclidean space with the matrix $\frac{\partial F}{\partial y}(x,y)$ only continuous at the base  point (thus, $\partial F/\partial y$ may be discontinuous elsewhere and $\partial F/\partial x$ may be everywhere discontinuous). In the case of a single scalar equation we show that this continuity hypothesis is unnecessary. The Inverse Function Theorem is shown as a consequence of the Implicit Function Theorem. Besides following Dini's inductive approach (see \cite{Dini}), these proofs avoid compactness arguments, fixed-point theorems, and integration theory. Instead of such tools, the proofs that follow employ the intermediate-value and the mean-value theorems, on the real line, and the intermediate-value property for derivatives on $\mathbb R$ (Darboux's property).

\vspace{0,2 cm}

Henceforth, we shall freely assume that all the functions are defined on a subset of a finite-dimensional Euclidean space. 

\vspace{0,2 cm}

Some comments are worthwhile concerning proofs of the implicit and inverse function theorems. Most proofs of the classical versions (enunciated for maps of class $C^1$ on an open set) start by showing the Inverse Function Theorem  and then derive the Implicit Function Theorem as a corollary. In general, these proofs employ either a compactness argument or the contraction mapping principle, 
see Krantz and Parks [10, pp. 41--52] and Dontchev and Rockafellar [4, pp. 9--20]. On the other hand, a proof of the classical Implicit Function Theorem that does not use either of these two tools can be seen in de Oliveira \cite{de Oliveira}.

\vspace{0,2 cm}

Taking into account maps that are everywhere differentiable (their differentials may be everywhere discontinuous), a proof of the Implicit Function Theorem can be found in Hurwicz and Richter \cite{Hurwicz}, whereas a proof of the Inverse Function Theorem can be seen in Saint Raymond \cite{Saint Raymond}. The first proof employs Brouwer's fixed-point theorem  while the second relies on Lebesgue's integration theory. Instead of assuming the continuity of the first order partial derivatives, these proofs assume an appropriate nondegeneracy condition at all points inside some open set containing the base point. It is worth noting that this 
quite general 
condition can be difficult to verify.

\vspace{0,2 cm}

Considering maps that are 
differentiable at the base point, but not necessarily differentiable elsewhere,
one can find proofs of the implicit and inverse function theorems in Hurwicz and Richter \cite{Hurwicz} and Nijenhuis \cite{Nijenhuis}. This last work extends Leach ~\cite{Leach} and employs the concept of {\sl strong differentiability} (also called {\sl strict differentiability}) and Banach's fixed-point theorem. It is worth noting that giving a differentiable map $F$ and a base point $p$, then the map $F$ is strong differentiable at $p$ if and only if its differential is continuous at $p$.

\vspace{0,2 cm}

Removing altogether the differentiability hypothesis, a proof of the Inverse Function Theorem for a Lipschitzian map can be seen in Clarke \cite{Clarke}. Yet, proofs of the Implicit Function Theorem for continuous maps can be found in Jittorntrum \cite{Jittorntrum} and Kumagai \cite{Kumagai}.

\vspace{0,1 cm}

In this article, the overall stategy of the proof of the Implicit Function Theorem is as follows. First, we prove it for a differentiable real function. Then, given a finite number of equations, we  prove it supposing that the matrix $\frac{\partial F}{\partial y}(x,y)$ is continuous at the base point. In addition, we prove the Inverse Function Theorem for a map whose Jacobian matrix is continuous at the base point.

\section{Notations and Preliminaries.}

Apart from the intermediate-value and the mean-value theorems, both on the real line, we assume the intermediate-value theorem for derivatives on $\mathbb R$ (Darboux's property): {\sl Given a differentiable function $f:[a,b]\to \mathbb R$, the image of the derivative function is an interval.}

\vspace{0,1 cm}

Let us consider $n$ and $m$, both in $\mathbb N$, 
and fix the canonical bases $\{e_1,\ldots ,e_n\}$ and $\{f_1, \ldots ,f_m\}$, of $\mathbb R^n$ and $\mathbb R^m$, respectively. Given $x=(x_1,\ldots, x_n)$ and $y=(y_1,\ldots,y_n)$, both in $\mathbb R^n$, we put $\left<x, y\right>= x_1y_1 + \cdots + x_ny_n$ and 
$|x|=\sqrt{\left<x,x\right>}$. Given $r>0$, let us write
$B(x;r)=\{y\ \textrm{in}\ \mathbb R^n:\ |y-x|<r\}$. 
We identify a linear map $T:\mathbb R^n \rightarrow \mathbb R^m$ with the $m\times n$ matrix  $M=(a_{ij})$, where $T(e_j)= a_{1j}f_1 + \cdots + a_{mj}f_m$ for  $j=1,\ldots,n$.
We also write $Tv$ for $T(v)$.

\vspace{0,1 cm}

In this section, $\Omega$ denotes a nonempty open subset of $\mathbb R^n$, where $n\geq 1$. Given a map $F: \Omega \rightarrow \mathbb R^m$ and a point $p$ in $\Omega$, we write $F(p)=\big(F_1(p),\ldots,F_m(p)\big)$. Let us suppose that $F$ is differentiable at $p$. 
The Jacobian matrix of $F$ at $p$ is 
\[JF(p)=\left(\frac{\partial F_i}{\partial x_j}(p)\right)_{\substack{1\leq i\leq m\\1 \leq j\leq n}}=
\left(\begin{array}{lllll}
\frac{\partial F_1}{\partial x_1}(p) & \cdots & \frac{\partial F_1}{\partial x_n}(p)\\
\ \ \ \vdots  &   & \ \ \ \vdots \ \ \      \\
\frac{\partial F_m}{\partial x_1}(p) & \cdots & \frac{\partial F_m}{\partial x_n}(p)
\end{array}
\right).\]
If $F$ is a real function, then we have $JF(p)= \nabla F(p)$, the gradient of $F$ at $p$.

\vspace{0,1 cm}

The following lemma (a particular case of the chain rule but sufficient for our purposes) is a local result. For practical reasons we state it for $\Omega=\mathbb R^n$. We omit the proof of the lemma.

\begin{lemma}\label{L1} Let $F:\mathbb R^n \to \mathbb R^m$ be  differentiable, $T: \mathbb R^k \to \mathbb R^n$ be the linear function associated to a $n\times k$ real matrix $M$, and $y$ be a fixed point in $\mathbb R^n$. Then, the function $G(x)= F(y + Tx)$, where $x$ is in $\mathbb R^k$, is differentiable and satisfies $JG(x)=JF(y + Tx)M$, for all $x$ in $\mathbb R^k$.
\end{lemma}
Given $a$ and $b$, both in $\mathbb R^n$, we put  $\overline{ab}=\{a+t(b-a):0\leq t\leq 1\}$.
The following mean-value theorem (in several variables) is a trivial consequence of the mean-value theorem on the real line and thus we omit the proof.

\begin{lemma}\label{L2} Let us consider a differentiable real function $F:\Omega \to \mathbb R$, with $\Omega$ open in $\mathbb R^n$. Let $a$ and $b$ be points in $\Omega$ such that the segment $\overline{ab}$ is within $\Omega$. Then, there exists $c$ in $\overline{ab}$ 
, with $c \neq a$ and $c\neq b$, that satisfies
\[F(b)-F(a) = \left<\nabla F(c), b-a\right>.\]
\end{lemma}

We denote the determinant of a real square matrix $M$ by $\det M$.

\begin{lemma} \label{L3} Let us consider a differentiable map $F: \Omega \to \mathbb R^n$, with $\Omega$ open within $\mathbb R^n$, and $p$ a point in $\Omega$ satisfying $\det JF(p)\neq 0$. Let us suppose that the real function $\det \big(\frac{\partial F_i}{\partial x_j}(\xi_{ij})\big)$ in the $n^2$ variables $\xi_{ij}$, with $1\leq i,j\leq n$ and $\xi_{ij}$ running in $\Omega$, is continuous at the point defined by $\xi_{ij}=p$, for all $1\leq i,j\leq n$. Then, the restriction of $F$ to some non-degenerate open ball $B(p;r)$ is injective.
\end{lemma}
{\bf Proof.} Since $ \det \big(\frac{\partial F_i}{\partial x_j}(p)\big)\neq 0$, the continuity hypothesis yields a $r>0$ such that $\det\big(\frac{\partial F_i}{\partial x_j}(\xi_{ij})\big)$ does not vanish, for all $\xi_{ij}$ in $B(p;r)$ and $1\leq i,j\leq n$. 
Now, let $a$ and $b$ be distinct in $B(p;r)$. By employing 
the mean-value theorem in several variables to each component $F_i$ of $F$, we find $c_i$ in the segment $\overline{ab}$, within  $B(p;r)$, such that $F_i(b)-F_i(a)=\left<\nabla F_i(c_i),b-a\right>$. Hence, 
\[ \left(\begin{array}{l}
F_1(b)- F_1(a)\\
\ \ \ \ \ \ \ \ \, \vdots  \\
F_n(b) - F_n(a) \\
\end{array}
\right)= 
\left(\begin{array}{ccc}
\frac{\partial F_1}{\partial x_1}(c_1) &\cdots & \frac{\partial F_1}{\partial x_n}(c_1)\\
\   \vdots & & \vdots \\
\frac{\partial F_n}{\partial x_1}(c_n) &\cdots & \frac{\partial F_n}{\partial x_n}(c_n)\\
\end{array}
\right)
\left(\begin{array}{l}
b_1- a_1\\
\ \ \ \ \vdots  \\
b_n -a_n \\
\end{array}
\right).  \]
Thus, since $\det\big(\frac{\partial F_i}{\partial x_j}(c_i)\big)\neq 0$ and $b-a\neq 0$, we conclude that $F(b)\neq F(a)$. $\qed$

\vspace{0,2 cm}

Given a real function $F:\Omega\to \mathbb R$, a short computation shows that the following definition of differentiability is equivalent to that which is most commonly employed. We say that $F$ is differentiable at $p$ in $\Omega$ if there are a ball $B(p;r)$ within $\Omega$, with $r>0$,  a $v$ in $\mathbb R^n$, and a vector-valued  map $E:B(0;r)\to \mathbb R^n$ satisfying
$$ \left\{\begin{array}{ll}
F(p +h)= F(p) + \left<v,h\right> + \left<E(h),h\right>,\ \textrm{for all}\ |h|<r,\\
\textrm{where}\ E(0)=0\ \textrm{and}\  E(h)\to 0\ \textrm{as}\ h\to 0.
\end{array}
\right.
$$

\section{The Implicit  Function Theorem.}

\

The first implicit function result we prove concerns one equation, several variables and a differentiable real function whose partial derivatives need not be continuous at any point. 
In its proof, we denote the variable in $\mathbb R^{n+1}=\mathbb R^n\times \mathbb R$ by $(x,y)$, where $x=(x_1,\ldots,x_n)$ is in $\mathbb R^n$ and $y$ is in $\mathbb R$. Given a nonempty subset $X$ of $\mathbb R^n$ and a nonempty subset $Y$ of $\mathbb R$, it is well-known that the set $X\times Y=\{(x,y): x \in X\ \textrm{and}\ y\in Y\}$ is open in $\mathbb R^n\times \mathbb R$ if and only if $X$ and $Y$ are open.

\vspace{0,2 cm}

In the next theorem, $\Omega$ denotes a nonempty open set within $\mathbb R^n\times \mathbb R$.

\begin{theorem}\label{TEO1} Let $F:\Omega \to \mathbb R$ be differentiable, with $\frac{\partial F}{\partial y}$ nowhere vanishing, and $(a,b)$ a point in $\Omega$ such that $F(a,b)=0$. Then, there exists an open set $X\times Y$, within $\Omega$ and containing the point $(a,b)$, that satisfies the following.
\begin{itemize}
\item[$\bullet$] There exists a unique function $g:X\to Y$ that satisfies $F\big(x,g(x)\big)=0$, for all $x$ in $X$.
 
\item[$\bullet$] We have $g(a)=b$. The function $g:X \to Y$ is differentiable and satisfies  
\[\frac{\partial g}{\partial x_j}(x)= - \frac{\frac{\partial F}{\partial x_j}(x,g(x))}{\frac{\partial F}{\partial y}(x,g(x))},\ \textrm{for all}\ x \ \textrm{in}\ X, \ \textrm{where}\ j=1,\ldots,n.\]
\end{itemize}
Moreover, if $\nabla F(x,y)$ is continuous at $(a,b)$ then $\nabla g(x)$ is continuous at $x=a$.
\end{theorem}
{\bf Proof.} By considering the function $F(x+a, \frac{y}{c}+b)$, with $c=\frac{\partial F}{\partial y}(a,b)$, we may assume that $(a,b)=(0,0)$ and $\frac{\partial F}{\partial y}(0,0)=1$. Next, we split the proof into three parts: existence and uniqueness, continuity at the origin, and differentiability.

\begin{itemize}
\item[$\diamond$]{\sf Existence and Uniqueness.} Let us choose a non-degenerate $(n+1)$-dimensional parallelepiped $X\times [-r,r]$, centered at $(0,0)$ and within $\Omega$, whose edges are parallel to the coordinate axes and $X$ is open. Then, the function $\varphi(y)=F(0,y)$, where $y$ runs over $[-r,r]$, is differentiable with $\varphi'$ nowhere vanishing and $\varphi'(0)=1$. Thus, by Darboux's property we have $\varphi'>0$ everywhere and we conclude that $\varphi$ is strictly increasing. Hence, by the continuity of $F$ and shrinking $X$ (if necessary) we may assume that 
$$F\Big|_{X\times \{-r\}}<0\ \ \ \textrm{and}\ \ \ F\Big|_{X\times\{r\}}>0.$$
As a consequence, fixing an arbitrary $x$ in $X$, the function 
$$\psi(y)=F(x,y),\ \textrm{where}\ y\in [-r,r],$$
satisfies $\psi(-r)<0<\psi(r)$. Hence, by the  mean-value theorem there exists a point $\eta$ in the open interval $Y=(-r,r)$ such that $\psi'(\eta)=\frac{\partial F}{\partial y}(x,\eta)>0$. Therefore, by Darboux's property we have $\psi'(y)>0$ at every $y$ in $Y$. Thus, $\psi$ is strictly increasing and the intermediate-value theorem yields the existence of a unique $y$, we write $y=g(x)$, in the open interval $Y$ such that $F(x,g(x))=0$.

\item[$\diamond$]{\sf Continuity at  the origin.} Let $\delta$ satisfy $ 0 < \delta < r$. From above, there exists an open set $\mathcal{X}$, contained in $X$ and containing $0$, such that $g(x)$ is in the interval $(-\delta ,\delta)$, for all $x$ in $\mathcal{X}$. Thus, $g$ is continuous at $x=0$.

\item[$\diamond$]{\sf Differentiability.} From the differentiability 
of the real function $F$ at $(0,0)$, and writing $\nabla F(0,0)=(v,1)\in \mathbb R^n\times \mathbb R$ for the gradient of $F$ at $(0,0)$, it follows that there are functions $E_1:\Omega\to \mathbb R^n$ and $E_2:\Omega \to \mathbb R$ satisfying
$$\left\{\begin{array}{ll}
F( h,k)= \left<v, h\right> + k + \left<E_1(h,k),h\right> + E_2(h,k)k,\\
\\
\textrm{where}\ \lim\limits_{(h,k)\to (0,0)}E_j(h,k)= 0 = E_j(0,0), \ \textrm{for}\ j=1,2.

\end{array}
\right.
$$

Hence, substituting 
[we already proved that $g(h)\xrightarrow{h\to 0}g(0)=0$]
$$\left\{\begin{array}{lll}
k=g(h),\\
E_j\big(h,g(h)\big)=\epsilon_j(h),\ 
\textrm{with}\ \lim\limits_{h\to 0}\epsilon_j(h)=\epsilon_j(0)=0\ \textrm{for}\ j=1,2,
\end{array}
\right.
$$
and noticing that we have $F\big(h,g(h)\big)=0$, for all possible $h$, we obtain
$$\left<v,h\right> + g(h) \ + \left<\epsilon_1(h),h\right> +  \epsilon_2(h)g(h)=0.\ \ \ \ \ \ \ \ \ \ \  $$
Thus, 
$$ [1 + \epsilon_2(h)]g(h) = - \left<v,h\right> - \left<\epsilon_1(h),h\right>. $$
If $|h|$ is small enough, then we have $1 + \epsilon_2(h)\neq 0$ and we may write

$$g(h)= \left<-v,h\right> +\left<\epsilon_3(h),h\right>,$$
where 
$$\epsilon_3(h)= \frac{\epsilon_2(h)}{1+\epsilon_2(h)}v -\frac{\epsilon_1(h)}{1+\epsilon_2(h)} \ \textrm{and}\ \lim_{h\to 0}\epsilon_3(h)=0.$$
Therefore, $g$ is differentiable at $0$ and $\nabla g(0)=-v$. 

Now, given any $a'$ in $X$, we put $b'=g(a')$. Then, $g: X \to Y$ solves the problem $F\big(x,h(x)\big)=0$, for all $x$ in $X$, with the condition $h(a')= b'$. From what we have just done it follows that $g$ is differentiable at $a'$. $\qed$
\end{itemize}

Next, we prove the implicit function theorem for a finite number of equations. 
Some notations are appropriate. We denote the variable in 
$\mathbb R^n\times \mathbb R^m = \mathbb R^{n+m}$ by $(x,y)$, where $x=(x_1,\ldots,x_n)$ is in $\mathbb R^n$ and $y=(y_1,\ldots,y_m)$ in $\mathbb R^m$. Given a nonempty subset $X$ of $\mathbb R^n$ and a nonempty subset $Y$ of $\mathbb R^m$, it is well-known that the set $X\times Y=\{(x,y): x \in X \ \textrm{and}\ y \in Y\}$ is open in $\mathbb R^n\times \mathbb R^m$ if and only if $X$ and $Y$ are open.
Given $\Omega$ an open subset of $\mathbb R^n\times \mathbb R^m$ and a  
differentiable map
$F:\Omega \to \mathbb R^m$ we write   $F=(F_1,\ldots,F_m)$, with $F_i$ the ith component of $F$ and $i=1,\ldots, m$, and 
\[\frac{\partial F}{\partial y}=\left(\frac{\partial F_i}{\partial y_j}\right)_{\substack{1\leq i\leq m\\1 \leq j\leq m}}=
\left(\begin{array}{lllll}
\frac{\partial F_1}{\partial y_1} & \cdots & \frac{\partial F_1}{\partial y_m}\\
\ \ \ \vdots  &   & \ \ \ \vdots \ \ \      \\
\frac{\partial F_m}{\partial y_1} & \cdots & \frac{\partial F_m}{\partial y_m}
\end{array}
\right).\]
Analogously, we define the matrix $\frac{\partial F}{\partial x}=\big(\frac{\partial F_i}{\partial x_k}\big)$, where $1\leq i \leq m$ and $1\leq k \leq n$.

\begin{theorem}\label{TEO2}{\bf (The Implicit Function Theorem).} {\sl Let $F:\Omega \to \mathbb R^m$ be differentiable, where $\Omega$ is an open set in $\mathbb R^n\times \mathbb R^m$. Let us suppose that $(a,b)$ is a point in $\Omega$ such that $F(a,b)=0$ and $\det\frac{\partial F}{\partial y}(a,b)\neq 0$, with $\frac{\partial F}{\partial y}(x,y)$ continuous at $(a,b)$. Then, there exists an open set $X\times Y$, within $\Omega$ and containing $(a,b)$, satisfying the following conditions.
\begin{itemize}
\item[$\bullet$] There exists a unique function $g:X\to Y$ that satisfies $F\big(x,g(x)\big)=0$, for all $x$ in $X$.
\item[$\bullet$] We have $ g(a)=b$. Moreover, the map $g:X \to Y$ is differentiable and 
\[Jg(x)  = - \left[\frac{\partial F}{\partial y }(x,g(x))\right]_{m\times m}^{-1}\left[\frac{\partial F}{\partial x}(x,g(x))\right]_{m\times n},\ \textrm{for all}\ x \ \textrm{in}\ X.\]
\end{itemize}
In addition, if $JF(x,y)$ is continuous at $(a,b)$ then $Jg(x)$ is continuous at $x=a$.}
\end{theorem}
{\bf Proof.} 
Let us consider the invertible matrix $\frac{\partial F}{\partial y}(a,b)=M$ and the associated bijective linear function $\mathcal{M}:\mathbb R^m \to \mathbb R^m$. By employing Lemma \ref{L1} we conclude that the map $G(x,z)=F[x,b + \mathcal{M}^{-1}(z-b)]$,  defined on a small enough neighborhood of $(a,b)$,  satisfy $\frac{\partial G}{\partial z}(a,b)=MM^{-1}$ and the condition $G(a,b)=0$. Therefore, replacing $F$ by $G$ if necessary, we may suppose without loss of generality that $M$ is the identity matrix of order $m$.

\vspace{0,2 cm}

Next, we split the proof into four parts: finding $Y$, existence and differentiability, differentiation formula, and uniqueness.

\begin{itemize}
\item[$\diamond$] {\sf Finding $Y$.} Defining $\Phi(x,y)=\big(x,F(x,y))$, where $(x,y)$ is in $\Omega$, we have
$$J\Phi(x,y)=\left(\begin{array}{l|l}
\ I  &\ 0\\ 
\hline
\frac{\partial F}{\partial x}  & \frac{\partial F}{\partial y}
\end{array} \right)\ \textrm{and}\ \det J\Phi(x,y)= \det \frac{\partial F}{\partial y}(x,y), $$ 
with $I$ the identity matrix of order $n$ and $0$ the $n\times m$ zero matrix. Thus, 
$\det J\Phi(a,b)\neq 0$. By hypothesis the matrix 
$\frac{\partial F}{\partial y}(x,y)$ is continuous at $(a,b)$.
Next, in order to apply Lemma \ref{L3} we introduce the variables  $\xi_{lk}$ in $\Omega$, where $l$ and $k$ run in 
$\{1,\ldots, m+n\}$, and the notation $(z_1,\ldots,z_n,z_{n+1},\ldots, z_{n+m})=(x_1,\ldots,x_n,y_1,\ldots,y_m)$. 
Then, the real function
$\det\big(\frac{\partial \Phi_{l}}{\partial z_k}(\xi_{lk})\big)=\det\big(\frac{\partial F_i}{\partial y_j}(\xi_{i+n,j+n})\big)$ is continuous at the point defined by $\xi_{lk}=(a,b)$, for all $l,k=1,\ldots, m+n$. Therefore, by Lemma \ref{L3}
 and shrinking $\Omega$ if necessary, we may assume that $\Phi$ is an injective map. 
We may also assume that $\Omega$ is an open non-degenerate 
parallelepiped $\mathcal{X}_1\times Y$ centered at $(a,b)$ whose edges are parallel to the coordinate axes. Hence, $\mathcal{X}_1$ and $Y$ are open (parallelepipeds).
\item[$\diamond$]{\sf Existence and differentiability.} We claim that the system
$$\left\{\begin{array}{ll}
F_1(x,y_1,\ldots,y_m) =0,\\
F_2(x,y_1,\ldots,y_m) =0,\\
\ \ \ \ \ \ \  \ \ \ \ \ \ \vdots\\
F_m(x,y_1,\ldots,y_m) =0,\\ 
\end{array}
\right. \  \textrm{with the conditions}\ \ 
\left\{\begin{array}{llll}
y_1(a)=b_1\\
y_2(a)=b_2\\
\ \ \ \   \ \  \vdots\\
y_m(a)=b_m,\\
\end{array}
\right.
$$
has a differentiable solution $g(x)=\big(g_1(x),\ldots, g_m(x)\big)$ on some open set $X$ containing $a$ [i.e., we have $F\big(x,g(x)\big)=0$ for all $x$ in $X$ and $g(a)=b$].

Let us prove it by induction on $m$. The case $m=1$ follows from Theorem \ref{TEO1} since $\frac{\partial F}{\partial y}(a,b)=1$ and, by continuity, we can assume  $\frac{\partial F}{\partial y}\neq 0$ everywhere.

Assuming that the claim holds for $m-1$, let us examine the case $m$. Then,
 given a pair $(x,y)=(x,y_1,\ldots, y_m)$ we introduce the helpful notations $y'=(y_2,\ldots,y_m)$, $y=(y_1,y')$, and $(x,y)=(x,y_1,y')$.

\vspace{0,2 cm}

Next, let us consider the equation $F_1(x,y_1,y')=0$, where $x$ and $y'$ are independent variables and $y_1$ is the dependent variable, with the condition $y_1(a,b')=b_1$. 
Since $\frac{\partial F_1}{\partial y_1}(a,b_1,b')=1$, by continuity we may assume that the function $\frac{\partial F_1}{\partial y_1}(x,y_1,y')$ does not vanish. Hence,
by Theorem 1 there exists a differentiable function $\varphi(x,y')$ on some open set [let us say, $\mathcal{X}_2\times \mathcal{Y}'$] containing $(a,b')$ that satisfies
\[F_1[x,\varphi(x,y'),y']=0\ \textrm{(on $\mathcal{X}_2\times \mathcal{Y}'$)} \  \textrm{and the condition}\ \varphi(a,b')=b_1.\]
As a consequence, $\varphi(x,y')$ also satisfies the $m-1$ equations
\[\frac{\partial F_1}{\partial y_1}[x,\varphi(x,y'),y']\frac{\partial \varphi}{\partial y_j}(x,y') + \frac{\partial F_1}{\partial y_j}[x,\varphi(x,y'),y'] =0,\ \textrm{for}\ j=2,\ldots,m.\]

Thus, since $\frac{\partial F_1}{\partial y}=\big(\frac{\partial F_1}{\partial y_1},\ldots,\frac{\partial F_1}{\partial y_m} \big)$ is continuous at $(a,b_1,b')$, with $\frac{\partial F_1}{\partial y_1}$ nowhere vanishing, and $\varphi$ is continuous, with $\varphi(a,b')=b_1$, we conclude that 
$\frac{\partial \varphi}{\partial y'}=\big(\frac{\partial \varphi}{\partial y_2}, \ldots, \frac{\partial \varphi}{\partial y_m}\big)$ is continuous at $(a,b')$.

Now, we look at solving the system with $m-1$ equations
 
$$\left\{\begin{array}{lll}
F_2[x,\varphi(x,y'),y']=0\\
\ \ \ \ \ \ \ \ \ \ \ \ \vdots\\
F_m[x,\varphi(x,y'),y']=0\\ 
\end{array}
\right. , 
\  \textrm{with the condition}\ \ 
y'(a)=b'.
$$
Let us define $\mathcal{F}_i(x,y')=F_i[x,\varphi(x,y'),y']$, with $i=2,\ldots,m$, and  write $\mathcal{F}=(\mathcal{F}_2,\ldots,\mathcal{F}_m)$.
Then, since the entries of the matrices $\frac{\partial \varphi}{\partial y'}(x,y')$ and $\frac{\partial F}{\partial y}(x,y)$ are continuous at $(a,b')$ and $(a,b)$, respectively, with $\varphi(a,b')=b_1$, we conclude that the entries of 
$\frac{\partial \mathcal{F}}{\partial y'}(x,y')$
are continuous at $(a,b')$. 
Yet, by hypothesis $\frac{\partial F}{\partial y}(a,b)$ is the identity matrix of order $m$ and thus we find
\[\frac{\partial \mathcal{F}_i}{\partial y_j}(a,b')= \frac{\partial F_i}{\partial y_1}(a,b)\frac{\partial \varphi}{\partial y_j}(a,b') + \frac{\partial F_i}{\partial y_j}(a,b)= 0 + \frac{\partial F_i}{\partial y_j}(a,b),\ \textrm{for}\ 2\leq i,j\leq m
.\]
This shows that  the matrix $\frac{\partial \mathcal{F}}{\partial y'}(a,b')$ is the identity one, of order $m-1$. Therefore, by induction hypothesis there exists a differentiable function $\psi$ on an open set $X$ containing $a$  [with  $\psi(X)$ contained in $\mathcal{Y}'$] that satisfies 
\[ \left\{\begin{array}{ll}
F_i[x,\varphi\big(x,\psi(x)\big),\psi(x)\big]=0, \ \textrm{for all}\ x \ \textrm{in}\ X, \ \textrm{for all} \ i=2,\ldots,m,\\ 
\textrm{and the condition} \ \psi(a)=b'.
\end{array}
\right.\] 
Clearly, we also have $F_1\big[x,\varphi\big(x,\psi(x)\big),\psi(x)\big]=0$, for all $x$ in $X$. Defining $g(x) = \big(\varphi(x,\psi(x)),\psi(x)\big)$, with $x$ in $X$, we obtain $F[x,g(x)]=0$, for all $x$ in $X$, and $g(a)=\big(\varphi(a,b'),b'\big)=(b_1,b')=b$, with $g$ differentiable on $X$.
\item[$\diamond$] {\sf Differentiation formula.} Differentiating $F[x,g(x)]=0$ we find
\[\frac{\partial F_i}{\partial x_k} + \sum_{j=1}^m\frac{\partial F_i}{\partial y_j}\frac{\partial g_j}{\partial x_k}=0,\ \textrm{with}\ 1\leq i\leq m\ \textrm{and}\ 1\leq k\leq n.\]
In matricial form, we write $\frac{\partial F}{\partial x}\big(x,g(x)\big) + \frac{\partial F}{\partial y}\big(x,g(x)\big)Jg(x)=0$. 

\item[$\diamond$] {\sf Uniqueness.} If $h:X\to Y$ 
and $x$ in $X$ satisfy $F(x,h(x))=0$, we find
$\Phi(x,h(x))=(x,0)=\Phi(x,g(x))$. In the first part of this proof (the ``finding $Y$'' part) we established that $\Phi$ is injective. Thus, $h(x)=g(x)$.
\end{itemize}
$\qed$

\section{The Inverse Function Theorem.}

\
\begin{theorem} \label{TEO3}{\bf (The Inverse Function Theorem).} Let $F:\Omega \to \mathbb R^n$ be differentiable, where $\Omega$ is an open set in $\mathbb R^n$. Let us suppose that $x_0$ is a point in $\Omega$ such that $JF(x_0)$ is invertible, with $JF(x)$ continuous at $x_0$. Then, there exist an open set $X$ containing $x_0$, an open set $Y$ containing $y_0=F(x_0)$, and a differentiable function $G:Y\to X$ that satisfies $F\big(G(y)\big)=y$, for all $y$ in $Y$, and $G\big(F(x)\big)=x$, for all $x$ in $X$. In addition,
\[JG(y)= JF\big(G(y)\big)^{-1}, \ \textrm{for all} \ y \ \textrm{in}\ Y,\]
and $JG(y)$ is continuous at $y=y_0$.
\end{theorem}
{\bf Proof.} 
By Lemma \ref{L3} we may assume that $F$ is injective. The map $\Phi(y,x)= F(x)-y$, where $(y,x)$ runs over $\mathbb R^n\times \Omega$, is differentiable and $\Phi(y_0,x_0)=0$. Yet, 
$\frac{\partial \Phi}{\partial x}\big(y_0,x_0\big)=JF(x_0)$ is invertible and $J\Phi(y,x)$ is continuous at $(y_0,x_0)$. The Implicit Function Theorem guarantees an open set $Y$ containing $y_0$ and a differentiable map $G:Y \to \Omega$, with $JG(y)$ continuous at $y=y_0$, satisfying  
$$F\big(G(y)\big)=y,\ \textrm{for all}\ y \ \textrm{in}\ Y.$$

Thus, $G$ is bijective from $Y$ to $X=G(Y)$ and $F$ is bijective from $X$ to $Y$. We also have $X=F^{-1}(Y)$. Since $F$ is continuous, $X$ is open (and contains $x_0$).

Putting $F(x)=\big(F_1(x),\ldots,F_n(x)\big)$ and $G(y)=\big(G_1(y),\ldots,G_n(y)\big)$ and differentiating  $\big(F_1(G(y)),\ldots,F_n(G(y))\big)$ we find
\[ \sum_{k=1}^n\frac{\partial F_i}{\partial x_k}\frac{\partial G_k}{\partial y_j}= \frac{\partial y_i}{\partial y_j}=
\left\{\begin{array}{ll}
1,\ \textrm{if} \ i=j,\\
0, \ \textrm{if}\ i\neq j.
\end{array}
\right.\]
 $\qed$

\paragraph{Acknowledgments.}  The author is greatly indebted to Professors Robert B. Burckel and James V. Ralston for their very valuable comments and suggestions.

\bigskip


\bigskip

\end{document}